\documentclass[12pt]{article}
\usepackage{amssymb}
\begin{document}
\def\dis{\displaystyle}
\newcommand{\rd}{\mbox{\rm Rad}}
\newcommand{\kn}{\mbox{\rm ker}}
\newcommand{\psp}{\vspace{0.4cm}}
\newcommand{\pse}{\vspace{0.2cm}}
\newcommand{\ptl}{\partial}
\newcommand{\dlt}{\delta}
\newcommand{\Dlt}{\Delta}
\newcommand{\sgm}{\sigma}
\newcommand{\al}{\alpha}
\newcommand{\be}{\beta}
\newcommand{\G}{\Gamma}
\newcommand{\gm}{\gamma}
\newcommand{\lmd}{\lambda}
\newcommand{\td}{\tilde}
\newcommand{\ad}{\mbox{\rm ad\,}}
\newcommand{\stl}{\stackrel}
\newcommand{\ol}{\overline}
\newcommand{\es}{\epsilon}
\newcommand{\la}{\langle}
\newcommand{\ra}{\rangle}
\newcommand{\vf}{\varphi}
\newcommand{\vsi}{\varsigma}
\newcommand{\ves}{\varepsilon}
\newcommand{\vt}{\vartheta}
\newcommand{\wt}{\mbox{\rm wt}\:}
\newcommand{\sym}{\mbox{\rm sym}}
\newcommand{\for}{\mbox{\rm for}}
\newcommand{\mbb}{\mathbb}
\def\qed{\hfill \hfill \ifhmode\unskip\nobreak\fi\ifmmode\ifinner\else\hskip5pt
\fi\fi
 \hbox{\hskip5pt\vrule width4pt height6pt depth1.5pt\hskip 1 pt}}
 \def\Box{\qed}
\def\der{\mbox{\rm Der\,}}
\def\a{\alpha}
\def\b{\beta}
\def\sF{\hbox{$\sc I\hskip -3.5pt F$}}
\def\Z{\hbox{$Z\hskip -5.2pt Z$}}
\def\sZ{\hbox{$\sc Z\hskip -4.2pt Z$}}
\def\Q{\hbox{$Q\hskip -5pt\vrule height 6pt depth 0pt\hskip 6pt$}}
\def\R{\hbox{$I\hskip -3pt R$}}
\def\C{\hbox{$C\hskip -5pt \vrule height 6pt depth 0pt \hskip 6pt$}}
\def\J {\vec J}
\def\d{\delta}
\def\D{\Delta}
\def\g{\gamma}
\def\G{\Gamma}
\def\l{\lambda}
\def\L{\Lambda}
\def\o{\omiga}
\def\p{\psi}
\def\Si{\Sigma}
\def\si{\sigma}
\def\sc{\scriptstyle}
\def\ssc{\scriptscriptstyle}
\def\dis{\displaystyle}
\def\cl{\centerline}
\def\nl{\newline}
\def\DD{{\cal D}}
\def\ll{\leftline}
\def\rl{\rightline}
\def\sF{\hbox{$\sc I\hskip -2.5pt F$}}
\def\ol{\overline}
\def\ul{\underline}
\def\wt{\widetilde}
\def\wh{\widehat}
\def\rar{\rightarrow}
\def\Rar{\Rightarrow}
\def\lar{\leftarrow}
\def\Lar{\Leftarrow}
\def\rla{\leftrightarrow}
\def\Rla{\Leftrightarrow}
\def\bs{\backslash}
\def\hs{\hspace*}
\def\vs{\vspace*}
\def\rb{\raisebox}
\def\ra{\rangle}
\def\la{\langle}
\def\Rad{\mbox{\rm Rad}}
\def\SS{\hbox{$S\hskip -6.2pt S$}}
\def\hi{\hangindent}
\def\ha{\hangafter}
\def\ni{\noindent}
\def\AA{{\cal A}}
\def\BB{{\cal B}}
\def\CC{{\cal C}}
\def\JJ{{\cal J}}
\def\KK{{\cal K}}

\def\v{\vec}
\def\SGN{{\rm sgn}}
\def\HOM{\mbox{\rm Hom}'_{\sZ}(\G,\mbb{F})}
\def\hom{\mbox{\rm Hom}^*_{\sZ}(\G,\mbb{F})}

\def\vi{\vec i}
\def\vj{\vec j}
\def\vk{\vec k}
\def\vm{\vec m}
\def\ii{{\bf i}}
\def\jj{{\bf j}}
\def\kk{{\bf k}}
\def\sone{{1\hskip -5.5pt 1}}
\def\one{{1\hskip -6.5pt 1}}
\def\sJ{J}
\def\th{\theta}
\def\isom{\mbox{${}^{\cong}_{\rar}$}}
\def\N{\mathbb{N}}
\def\Z{\mathbb{Z}}
\def\sZ{\mathbb{Z}}
\def\Q{\mathbb{Q}}
\def\R{\mathbb{R}}
\def\C{\mathbb{C}}
\def\sC{\mathbb{C}}
\def\sF{\mathbb{F}}
\def\F{\mathbb{F}}
\cl{{\bf\large Derivations and 2-Cocycles of Contact Lie
Algebras}}\cl{{\bf\large Related to Locally-Finite
Derivations}\footnote{\,1991 Mathematical Subject Classification.
Primary 17B65; Secondary 58F05\nl\hspace*{4ex}\,Supported by NSF
grant 10171064 of China and two grants ``Excellent Young Teacher
Program'' and ``Trans-Century Training Programme Foundation for the
Talents'' from Ministry of Education of China}} \vskip3pt
\cl{(appeared in {\it Comm. Alg.} {\bf32} (2004), 4613--4631.)}
\par
\vs{3pt}\par \cl{Guang'an Song and Yucai
Su}\vs{1pt}\par\cl{\small Department of Mathematics, Shanghai
Jiaotong University}\cl{\small Shanghai
200030, P.~R.~China }%
\cl{\small Email: gasong@sjtu.edu.cn, ycsu@sjtu.edu.cn}%
{\vs{4pt}
\par\ni{\small%
{\bf Abstract}. Classical contact Lie algebras are the fundamental
algebraic structures on the manifolds of contact elements of
configuration spaces in classical mechanics. Xu introduced a large
category of contact simple Lie algebras which are related to
locally finite derivations and are in general not finitely graded.
The isomorphism classes of these Lie algebras were determined in a
previous paper by Xu and Su. In this paper, the derivation
algebras and the 2-cohomology groups of these Lie algebras are
determined, and it is obtained that the 2-cohomology groups of
these Lie algebras are all trivial.}
\vs{2pt}

\ni{\bf Key words:} Lie algebras of contact type, structure, derivation, $2$-cocycle
}
\par\
\vs{-5pt}\par
\cl{\bf\S1. Introduction}
\par
One of the four well-known classes of infinite dimensional simple Lie
algebras of Cartan type is the class of contact Lie algebras.
Classical contact Lie algebras are the fundamental algebraic structures on
the manifolds of contact elements of configuration spaces in classical
mechanics.
Simple contact Lie algebras have been studied by
Kac [2, 3], Osborn [5], Osborn and Zhao [6], and Xu [9, 10].
In [9], Xu constructed a class of contact Lie algebras based on the
pairs $(\AA,\DD)$, where $\AA$ is a commutative associative unital
algebra and $\DD$ is a finite dimensional commutative subalgebra
of locally finite derivations of $\AA$ such that $\AA$ does not have
proper $\DD$-stable ideal.
{}From the classification of such pairs $(\AA,\DD)$
in [8], it is known that the class of contact Lie algebras in [9] is the
largest class under this locally finite condition.
\par
In [7], Xu and Su determined the isomorphism classes of
the contact Lie algebras given in [9]. Osborn and Zhao [6]
determined the isomorphism classes of the contact Lie algebras
constructed by themselves. The problem of determining the
derivation algebras and the 2-cohomology groups of the Lie
algebras given in [6] remains unsolved. In this paper, we shall
solve this problem for the contact Lie algebras given in [9] (thus
also solve the problem for the Lie algebras given in [6]).
\par
The significance of derivations for Lie theory, as pointed in [1],
resides in their affinity to the cohomology groups, their
determination affords insight into structural features of Lie
algebras which do not figure prominently in the defining
properties. Some general results concerning derivations of
finitely graded Lie algebras were established in [1]. However in
our case the Lie algebras are in general not finitely graded, the
results of [1] can not apply here, thus the determination of the
derivation algebras of Lie algebras of this kind is a nontrivial
problem. We use a different technique to solve the problem here.
\par
The 2-cohomology groups of Lie algebras play important roles in
the central extensions of Lie algebras, which are often used in
the structure theory and the representation theory of Kac-Moody
algebras (cf.~[4]). Since the cohomology groups are closely
related to the structures of Lie algebras, the computation of
cohomology groups seems to be important and interesting as well.
\par
In Section 2, we shall present the contact simple Lie algebras
given in [9] in what we call the {\it normalized forms} (cf.~[7]).
Then the determinations of derivation algebras and 2-cohomology
groups of these Lie algebras are given in Sections 3 and 4 (see
Theorems 3.1 and 4.6) and we obtain that the 2-cohomology groups
of these Lie algebras are all trivial.
\par\
\vs{-2pt}\par
\cl{\bf\S2. Normalized forms}
\vs{-1pt}\par
For $m,n\in\Z$, we denote $\ol{m,n}=\{m,m+1,...,n\}$.
Let
$\vec\ell=(\ell_1,...,\ell_6)\in \mbb{N}{\ssc\,}^6$
such that $\sum_{p=1}^6\ell_p>0.$
We \vs{-3pt}denote
$$
\iota_0=0,\;\;\;\iota_i=\ell_1+\ell_2+...+\ell_i,\;\;\;\;
I_i=\ol{\iota_{i-1}+1,\iota_i}\;\;\;\;\for\;\;\;i=1,2,...,6.
\vs{-3pt}\eqno(2.1)$$
Set
$I_{i,j}=\cup_{i\le p\le j}I_p=\ol{\iota_{i-1}+1,\iota_j}$
for $1\le i\le j\le6.$
Set
$$
I=I_{1,6}=\ol{1,\iota_6},\;\;\;\;J=\ol{1,2\iota_6},\;\;\;\;
\wh K=\{0\}\cup K\;\;\;\;\for\;\;\;K\subset J.
\eqno(2.2)$$
Define the map \rb{5pt}{$\ol{\ }$}\,$:J\rar J$ to be the index
shifting of $\iota_6$ steps in $J$, i.e.,
$$
\ol p=p+\iota_6\mbox{ \ if \ }p\in\ol{1,\iota_6},
\mbox{ \ \ \,or \ }\ol p=p-\iota_6\mbox{ \ if \ }p\in\ol{\iota_6+1,2\iota_6}.
\eqno(2.3)$$
For any subset $K$ of $J$, we denote $\ol K=\{\ol p\,|\,p\in K\}.$
Thus $J=I\cup\ol I$. Set
$J_i=I_i\cup\ol I_i,$ $J_{i,j}=I_{i,j}\cup\ol I_{i,j}.$
\par
An element of the vector space $\mbb{F}^{1+2\iota_6}$ is denoted as
$$
\a=(\a_0,\a_1,\a_{\ol1},...,\a_{\iota_6},\a_{\ol\iota_6})
\mbox{ \ with \ }\a_p\in\mbb{F}\ \ \mbox{for all}\ \ p\in\wh J.
\eqno(2.4)$$
For $\a\in\mbb{F}^{1+2\iota_6}$ and $K\subset J$, we denote by $\a_{\ssc K}$
the vector obtained from $\a$ with support $K$, i.e.,
$$
\begin{array}{ll}
\a_{\ssc K}=
\!\!\!\!&
(\b_0,\b_1,\b_{\ol1},...,\b_{\iota_6},\b_{\ol\iota_6})
\vs{4pt}\\&
\mbox{ with }\
\b_p=0\mbox{ \ if \ }p\notin K\mbox{ \ and \ }\b_p=\a_p\mbox{ \ if \ }p\in K.
\end{array}
\eqno(2.5)$$
When the context is clear, we also use $\a_K$ to denote the element in
$\mbb{F}^{|K|}$ obtained from $\a$ by deleting the coordinates at
$p\in\ol{0,2\iota}\bs K$; for instance, $\a_{\{1,\ol2\}}=(\a_1,\a_{\ol2}).$
Moreover, we denote
$$
a_{[p]}=(0,...,0,\stl{p}{a},0,...,0)\in\mbb{F}^{1+2\iota}
\qquad\for\;\;a\in\mbb{F}.
\eqno(2.6)$$
\par
Take
$$
\si_{\ol p}=\si_p=\left\{\begin{array}{ll}
-1_{[p]}-1_{[\ol p]}&\mbox{if \ }p\in I_{1,3},\vs{2pt}\\
-1_{[p]}&\mbox{if \ }p\in I_{4,5},\vs{2pt}\\
0&\mbox{if \ }p\in\wh I_6.
\end{array}\right.
\eqno(2.7)$$
Let $\G$ be an additive subgroup of $\mbb{F}^{1+2\iota_6}$ such that
$$
\begin{array}{c}
\{1_{[p]}\,|\,p\in J_{1,3}\cup I_{4,5}\}\subset
\G\subset\{\a\in\mbb{F}^{1+2\iota_6}\,|\,\a_{I_6\cup\ol
I_{4,6}}=0\},\mbox{ and \ }\vs{4pt}\\%
1_{[0]}\in\G
 \mbox{ if }\G_0\ne\{0\},
\end{array}
\eqno\begin{array}{r}\vs{4pt}(2.8)\!\!\\ (2.9)\!\!\end{array}$$
where in general, we define $\G_p=\{\a_p\mid
(\a_0,\a_1,\a_{\ol1},...,\a_{\iota_6},\a_{\ol\iota_6})\in\G\}.$
\par
Take
$\JJ_0=\{0\}\mbox{ or }\mbb{N}$
such that $\JJ_0+\G_0\ne\{0\},$ and
$$\JJ_1=\{\ii\in\mbb{N}\:^{2\iota_6}\,|\,
\ii_{I_{1,2}\cup I_4\cup \ol I_1}=0\},
\eqno(2.10)$$
and set
$\JJ=(\JJ_0,\JJ_1)\subset \mbb{N}^{1+2\ell_6},$
an additive subsemigroup. An element of $\JJ$  is denoted as
$$
\vec i=(i_0,\ii)=(i_0)_{[0]}+\ii,\mbox{ with }i_0\in\JJ_0,\ii\in\JJ_1.
\eqno(2.11)$$
\par
Let $\AA$ be the semigroup algebra $\mbb{F}[\G\times\JJ]$ with
basis $B$ and multiplication $\cdot$ defined by
$$
\begin{array}{ll}
B=\{x^{\a,\vec i}\,|\,(\a,\vec i)\in\G\times \JJ\},\,\,\;
\vs{4pt}\\
 x^{\a,\vec i}\cdot x^{\b,\vec j}=x^{\a+\b,\vec i+\vec j}
\ \ \for\ \ (\a,\vec i),(\b,\vec j)\in\G\times \JJ.
\end{array}
\eqno(2.12)$$
Then $(\AA,\cdot)$ forms a commutative associative algebra with
$1=x^{0,0}$ as the identity element. For convenience, we often
denote
$$x^{\a}=x^{\a,0},\ \
t^{\vec i}=x^{0,\vec i},\ \
t_p=t^{1_{[p]}}\qquad \for\ p\in\wh J.
\eqno(2.13)$$
In particular,
$t^{\vec i}=\prod_{p\in\wh J}t_p^{i_p}.$
\par
Define the operators $\{\ptl_p,\ptl^*_p,\ptl_{t_p}\,|\,p\in\wh
J\}$ on $\AA$ by (in fact, they are derivations of $(\AA,\cdot)$)
$$\ptl_p=\ptl^*_p+\ptl_{t_p}\mbox{ \ and \ }
\ptl^*_p(x^{\a,\vec i})=\a_p x^{\a,\vec i},\ \
\ptl_{t_p}(x^{\a,\vec i})=i_p x^{\a,\vec i-1_{[p]}}, \eqno(2.14)$$
for $p\in\wh J,(\a,\vec i)\in\G\times\JJ,$ where we always use the
convention that $x^{\a,\vec i}=0$ if $(\a,\vec
i)\notin\G\times\JJ$. In particular,
$$\ptl^*_p=0,\ \ \ptl_{t_q}=0\ \ \for\ \ p\in
I_6\cup\ol I_{4,6},\
q\in I_{1,2}\cup I_4\cup \ol I_1,
\eqno(2.15)$$
by (2.8), (2.10). We call $\ptl^*_p$ a {\it grading operator} if
$\ptl^*_p\ne0$; $\ptl_{t_q}$ a {\it down-grading operators} if
$\ptl_{t_q}\ne0$; and $\ptl_r^*+\ptl_{t_r}$ a {\it mixed operator}
if $\ptl_r^*\ne0\ne\ptl_{t_r}$. Then the types of derivation pairs
in the order of the groups $\{(\ptl_p,\ptl_{\bar{p}})\,|\,p\in
I_i\}$ for $i\in\ol{1,6}$ are
$$
(g,g),\;(g,m),\;(m,m),\;(g,d),\;(m,d),\;(d,d),
\eqno(2.16)$$
where ``m'', ``g'', ``d'' stand for mixed, grading, down-grading operators.
\par
We denote
$$
\ptl=\sum_{p\in J_{1,3}\cup I_{4,5}} \ptl_p^*+\sum_{p\in
I_6\cup\ol I_{4,6}}t_p\ptl_{t_p}. \eqno(2.17)$$ Define the
following Lie bracket $[\cdot,\cdot]$ on $\AA$:
$$
\begin{array}{ll}
[u,v]=\!\!\!\!&\dis
 \sum_{p\in I}x^{\si_p} (\ptl_p(u)\ptl_{\ol p}(v)-\ptl_{\ol
p}(u)\ptl_p(v))
\vs{4pt}\\&
+ (2-\ptl)(u)\ptl_0(v)-\ptl_0(u)(2-\ptl)(v),
\end{array}
\eqno(2.18)$$
 for $u,v\in\AA$ (cf.~(2.3)). Then
$(\AA,[\cdot,\cdot])$ forms a contact Lie algebra, which is in
general not finitely-graded. The algebras $(\AA,[\cdot,\cdot])$
are the normalized forms of the contact simple Lie algebras given
in [9]. We denote the Lie algebra $(\AA,[\cdot,\cdot])$ by
$$
\KK=\KK(\vec\ell,\si,\G,\JJ),\mbox{ where } \si=\sum_{p\in
I_{1,5}}\si_p. \eqno(2.19)$$ Then the Lie algebras given in [6]
are the Lie algebras $\KK(\vec\ell,\si,\G,\{0\})$ with
$\vec\ell=(\ell_1,0,0,0,0,0)$.
%
\par\ \vs{-2pt}\par%
\cl{\bf\S3. Structure of derivation algebras}
\vs{-0.5pt}\par
In this section, we shall determine $\der\KK$.
Recall that a {\it derivation} $d$ of the Lie algebra $\KK$ is a linear
transformation on $\KK$ such that
$$
d([u_1,u_2])=[d(u_1),u_2]+[u_1,d(u_2)]\ \ \for\ \ u_1,u_2\in\KK.
\eqno(3.1)$$ Denote by $\der\KK$ the space of the derivations of
$\KK$, which is a Lie algebra. Moreover, the space $\ad\KK=\{\ad
u\,|\,u\in\KK\}$ is an ideal. Elements in $\ad\KK$ are called {\it
inner derivations}, while elements in $\der\KK\bs\ad\KK$ are
called {\it outer derivations}.
\par
We give the explicit form of (2.18) as follows:
$$
\begin{array}{ll}
[x^{\a,\dis\vec i},x^{\b,\vec j}]= \!\!\!\!&\dis \sum_{p\in
I_{1,3}}(\a_p\b_{\ol p}-\a_{\ol p}\b_p) x^{\si_p+\a+\b,\vec i+\vec
j}
\vs{4pt}\\&\dis+\sum_{p\in I_{2,5}}(\a_pj_{\ol p}-i_{\ol p}\b_p)
x^{\si_p+\a+\b,\vec i+\vec j-1_{[\ol p]}} \vs{4pt}\\ &\dis
+\sum_{p\in I_3}(i_p\b_{\ol p}-j_p\a_{\ol p})
x^{\si_p+\a+\b,\vec i+\vec j-1_{[p]}}
\vs{4pt}\\&\dis+\sum_{p\in I_3\cup
I_{5,6}}(i_pj_{\ol p}-i_{\ol p}j_p) x^{\si_p+\a+\b,\vec i+\vec
j-1_{[p]}-1_{[\ol p]}} \vs{4pt}\\ &\dis +((2-\vt(\a,\vec
i))\b_0-\a_0(2-\vt(\b,\vec j))) x^{\a+\b,\vec i+\vec j} \vs{4pt}\\
&\dis +((2-\vt(\a,\vec i))j_0-i_0(2-\vt(\b,\vec j))) x^{\a+\b,\vec
i+\vec j-1_{[0]}},
\end{array}
\eqno(3.2)$$
for $(\a,\vec i),(\b,\vec j)\in\G\times\J$, where
$$
\vt(\a,\vec i)=\sum_{p\in J_{1,3}\cup I_{4,5}}\a_p+
\sum_{p\in I_6\cup\ol I_{4,6}}i_p\qquad\for\ (\a,\vec i)\in\G\times\JJ.
\eqno(3.3)$$
In particular,
$$
[1,x^{\b, \vec j}]=2\b_0x^{\b,\vec j} +2j_0x^{\b,\vec j-1_{[0]}},
\eqno(3.4)$$
$$
[x^{-\si_p},x^{\b,\vec j}]=\left\{
\begin{array}{ll}
(\b_{\ol p}-\b_p)x^{\b,\vec j}& \for\ p\in I_1,
\vs{4pt}\\
(\b_{\ol p}-\b_p)x^{\b,\vec j}+j_{\ol p}x^{\b,\vec j-1_{[\ol p]}}
& \for\ p\in I_2,
\vs{4pt}\\
(\b_{\ol p}-\b_p)x^{\b,\vec j}-j_px^{\b,\vec j-1_{[p]}}
+j_{\ol p}x^{\b,\vec j-1_{[\ol p]}}& \for\ p\in I_3,
\end{array}
\right. \eqno(3.5)$$
$$
[x^{-\si_q,1_{[\ol q]}},x^{\b,\vec j}]=\left\{
\begin{array}{ll}
(-\b_q+j_{\ol q})x^{\b,\vec j}& \for\ q\in I_4,
\vs{4pt}\\
(-\b_q+j_{\ol q})x^{\b,\vec j}
-j_qx^{\b,\vec j-1_{[q]}}
& \for\ q\in I_5,
\end{array}
\right. \eqno(3.6)$$
$$
[t^{1_{[r]}+1_{[\ol r]}},x^{\b,\vec j}]=
(j_{\ol r}-j_r)x^{\b,\vec j}\ \for\ r\in I_6.
\eqno(3.7)$$
\par
We shall find out all possible derivations of $\KK$. First,
observe from definition (2.14) and (2.18) that for $p\in\ol
I_2\cup J_3\cup I_5$, the operator $\ptl_{t_p}$ is an outer
derivation of $\KK$. Next, denote by $\HOM$ the set of group
homomorphisms $\mu:\G\to\mbb{F}$ such that $\mu(\si_p)=0$ for
$p\in I_{1,5}$. For $\mu\in\HOM$, we define a linear
transformation $d_\mu$ on $\KK$ by
$$
d_{\mu}(x^{\a,\v i})=\mu(\a)x^{\a,\v i}\;\;\;\for\;\;\;(\a,\v i)
\in\G\times\JJ. \eqno(3.8)$$ Clearly, by (3.2), $d_\mu$ is a
derivation of $\KK$. We identify $\HOM$ with a subspace of
$\der\,\KK$ by $\mu\mapsto d_\mu$. For $p\in I_{1,3}$, we define
$\mu_p\in\HOM$ by $\mu_p(\a)=\a_{\ol p}-\a_p$. By (3.5), we have
$$
d_{\mu_p}=\ad{x^{-\si_p}}+\ptl_{t_p}-\ptl_{t_{\ol p}}\;\;\;\;
\for\;\;\;p\in I_{1,3}.%
\eqno(3.9)$$%
We also define $\mu'_0:\a\mapsto\a_0$. If $\JJ_0=\N$, we see that
$d_{\mu'_0}={1\over2}(\ad 1-\ptl_{t_0})$ is an outer derivation;
in this case we set $\mu_0=0$. If $\JJ_0=\{0\}$, then
$d_{\mu'_0}=\ad 1$ is an inner derivation; in this case, we set
$\mu_0=\mu'_0$.
 We
fix a subspace $\hom$ of $\HOM$ such that
$$
\HOM=\hom\oplus{\rm span}\{\mu_p\,|\,p\in\wh I_{1,3}\},
\eqno(3.10)$$ is a direct sum as vector spaces. Then we have
\par\ni
{\bf Theorem 3.1}. {\it The derivation algebra $\der\,\KK$ is a
direct of subspaces:
$$
\der\,\KK=\bigoplus_{p\in\ol I_2\cup J_3\cup I_5}\mbb{F}\ptl_{t_p}
\bigoplus\hom\bigoplus\ad\KK.
\eqno(3.11)$$}
{\it Proof.}
First we prove that the right-hand side of (3.11) is the direct
sum. Thus suppose
$$
d=\sum_{p\in\ol I_2\cup J_3\cup
I_5}c_p\ptl_{t_p}+d_\mu+\sum_{(\a,\v i)\in\G\times\JJ}c_{\a,\v
i\,\,}\ad x^{\a,\v i},
$$
is the zero derivation, where $c_p,c_{\a,\v i}\in\F,\mu\in\hom$
such that $\{(\a,\v i)\in\G\times\JJ\,|\,c_{\a,\v i}\ne0\}$ is a
finite set. Applying $d$ to $x^{-\si_q,1_{[\ol q]}},t^{2_{[\ol
q]}}$ for $q\in
I_{4,5}$, we obtain that%
$$c_{\a,\v i}\ne0\ \Rar\ \a_q=i_{\ol q}=0\;\;\;\;\for\;\;\;\;q\in I_{4,5}.%
$$%
In particular, $c_{-\si_q,1_{[\ol q]}}=0$ for $q\in I_{4,5}$.
Applying $d$ to $x^\b$ for $\b\in\G$, by calculating the
coefficients of $x^{\b,\v j}$ with $\sum_{p\in\wh J}j_p$ being
maximal and by calculating the coefficient of $x^\b$, we obtain
$$\mu(\b)+c_{0,0}\b_0+\sum_{p\in I_{1,3}}c_{-\si_p,0}(\b_p-\b_{\ol
p})=0\;\;\;\;\mbox{for \ all}\;\;\;\;\b\in\G.%
\eqno(3.12)$$%
If $\JJ_0=\{0\}$, then by (3.10), we obtain
$\mu=0,\,c_{0,0}=c_{-\si_p}=0$ for $p\in I_{1,3}$; if $\JJ_0=\N$,
applying $d$ to $t_0$, we obtain in particular $c_{0,0}=0$, then
(3.10) and (3.12) again give $\mu=0,\,c_{-\si_p}=0$ for $p\in
I_{1,3}$. Now applying $d$ to the
set%
$$
\begin{array}{ll}
\{1,s,x^{-\si_p},x^{2_{[p]}},x^{2_{[\ol p]}},x^{-\si_q,1_{[\ol
q]}},x^{2_{[q]}},t^{2_{[\ol q]}},t^{1_{[r]}+1_{[\ol
r]}},
\!\!\!\!&
t^{2_{[r]}},t^{2_{[\ol r]}}\,|
\vs{4pt}\\&
\,p\in I_{1,3},q\in I_{4,5},r\in I_6\},%
\end{array}
$$%
where $s=0$ if $\G_0\ne\{0\}$ and $s=t_0$ if $\G_0=\{0\}$, we
obtain by induction on $\sum_{p\in\wh J}i_p$ that $c_{\a,\v i}=0$
for all $(\a,\v i)\in\G\times\JJ$. Finally applying $d$ to $t_p$
for $p\in\ol I_2\cup J_3\cup I_5$, we obtain $c_p=0$. Thus the
right-hand side of (3.11) is a direct sum.
\par
Now let $d\in\der\,\KK$ and let $D$ be the space in the right-hand
side of (3.11). Note that $D\supset\HOM$ by (3.9) and (3.10). We
shall prove that after a number of steps in each of which $d$ is
replaced by $d-d'$ for some $d'\in D$ the zero derivation is
obtained and thus proving that $d\in D$. This will be done by a
number of claims.
\par
{\bf Claim 1}. Let $$ A'_1=\{x^{-\si_p},x^{-\si_q,1_{[\ol
q]}}\,|\,p\in I_{2,3},q\in I_5\},\eqno(3.13)$$ and $A_1=A'_1$ if
$\JJ_0=\{0\}$ or $A_1=A'_1\cup\{1\}$ if $\JJ_0=\N$ (note that
$A_1$ is the set of {\it ad}-locally finite but not {\it
ad}-semisimple elements of the basis $B$). Say $A_1$ has $n_1$
elements, and list its elements as $y_1,...,y_{n_1}$. Let
$u\in\KK$. For any $r\in\ol{1,n_1}$, we can choose $v\in\KK$
satisfying $u=[y_r,v]$ such that if $[y_q,u]=0$, then $[y_q,v]=0$
for $q\in\ol{1,r-1}$, i.e., any element $u$ of $\KK$ is an image
of the operator $\ad y_r$ such that if $u$ commutes with $y_q$,
then we can find a preimage $v$ also commuting with $y_q$.
\par
Since $u$ is a linear combination of the basis elements $x^{\a,\v
i}$, we can suppose $u=x^{\a,\v i}$. We shall use (3.4)-(3.6) to
prove the claim. Say $y_r=x^{-\si_p}$ for some $p\in I_3$
(otherwise the proof is
similar).%
\par%
If $\a_{\ol p}\ne\a_p$, by (3.5), we can take $v_1=(\a_{\ol
p}-\a_p)^{-1}x^{\a,\v i}$ such that
$$u-[y_r,v_1]=u_1,\mbox{ where }
u_1=(\a_{\ol p}-\a_p)^{-1} (i_px^{\a,\v i-1_{[p]}}-i_{\ol
p}x^{\a,\v i-1_{[\ol p]}}).%
$$%
By induction on $i_p+i_{\ol p}$, we can take $v_2$ such that
$u_1=[y_r,v_2]$. Thus we can find $v=v_1+v_2$ such that
$u=[y_r,v]$. Furthermore, if $[y_q,u]=0$, we must also have
$[y_q,v_1]=0$ (by (3.4)-(3.6), since $y_q\in A_1$), and by
induction, we also have $[y_q,v_2]=0$.%
\par%
Now suppose $\a_{\ol p}=\a_p$. Take $v_1=(i_{\ol
p}+1)^{-1}x^{\a,\v i+1_{[\ol p]}}$, we have
$$
u-[y_r,v_1]=v_2,\mbox{ where }v_2=(i_{\ol p}+1)^{-1}i_px^{\a,\v
i+1_{[\ol p}]-1_{[p]}}.%
$$%
By induction on $i_p$, we also have the result. This proves the
claim.
\par
{\bf Claim 2}.  We can suppose $d(A_1)=0$.
\par
For $r\in\ol{1,n_1}$, suppose we have proved $d(y_q)=0$ for $q<
r$. Let $u=d(y_r)$. We have $[y_q,u]=d([y_q,y_r])-[d(y_q),y_r]=0$
(note that $A_1$ is commutative). By Claim 1, there is $v\in\KK$
such that $u=[y_r,v]$ and $[y_q,v]=0$. Replacing $d$ by $d+\ad v$,
we obtain $d(y_q)=0$ for $q\le r$, thus the claim follows.
\par
{\bf Claim 3}. Let $$A'_2=\{x^{-\si_p},x^{-\si_q,1_{[\ol
q]}},t^{1_{[r]}+1_{[\ol r]}}\,|\,p\in I_1,q\in I_4,r\in I_6\},
\eqno(3.14)$$ and let $A_2=A'_2\cup\{1\}$ if $\JJ_0=\{0\}$ or
 $A_2=A'_2$ if $\JJ_0=\N$ (note that elements of $A_2$ are all {\it ad}-semisimple
 and that $A_1\cup A_2$ is commutative). We can suppose $d(A_2)=0$.
\par
Say $A_2$ has $n_2$ elements and list its elements as
$z_1,...,z_{n_2}$. For $r\in\ol{1,n_2}$, suppose we have proved
$d(z_q)=0$ for $q< r$. Let $u=d(z_r)$. We have
$$
[y_p,u]=d([y_p,z_r])-[d(y_p),z_r]=0\mbox{ \ for \
}p\in\ol{1,n_1},%
\vs{-6pt} \eqno(3.15)$$
$$
[z_q,u]=d([z_q,z_r])-[d(z_q),z_r]=-[d(z_q),z_r]=0\mbox{ \ for \
}q<r.%
\eqno(3.16)$$%
 Write $u$ as a linear combination of the basis
$B$ (cf.~(2.12)). Suppose a term $x^{\a,\v i}$ appears in $u$ with
coefficient $c_{\a, \v i}\ne0$. The above two equations give
$$
[y_p,x^{\a,\v i}]=[z_q,x^{\a,\v i}]=0.\eqno(3.17)
$$
Say $z_r=x^{-\si_{r_1},1_{[\ol {r}_{_{\ssc1}}]}}$ for some $r_1\in
I_4$ (otherwise the proof is similar). If $i_{\ol
r_1}\ne\a_{r_1}$, then by replacing $d$ by $d+\ad v$ for
$v=c_{\a,\v i}(i_{\ol r_1}-\a_{r_1})^{-1}x^{\a,\v i}$, we see from
(3.6) that the term $x^{\a,\v i}$ then disappears in $d(z_r)$, and
furthermore, after this replacement, we see from (3.17) that
$d(y_p)=d(z_q)=0$ still holds. Thus we can suppose
$$u=d(z_r)\in{\rm span}\{x^{\a,\v i}\,|\,i_{\ol
r_1}=\a_{r_1}\}.\eqno(3.18)$$ Now for any $x^{\b,\v j}\in\KK$ with
$\b_{r_1}=j_{\ol r_1}$, we have
$$
\begin{array}{ll}
[x^{\b,\v j},u]
\!\!\!\!&
=d([x^{\b,\v j},x^{-\si_{r_1},1_{[\ol
r_{_{\ssc1}}]}}])-[d(x^{\b,\v j}),x^{-\si_{r_1},1_{[\ol
r_{_{\ssc1}}]}}]
\vs{4pt}\\&
=-[d(x^{\b,\v j}),x^{-\si_{r_1},1_{[\ol
r_{_{\ssc1}}]}}],\vs{-6pt}%
\end{array}
\eqno(3.19)$$
$$
\begin{array}{ll}
[x^{2_{[r_1]}},u]
\!\!\!\!&
=d([x^{2_{[r_1]}},x^{-\si_{r_1},1_{[\ol
r_{_{\ssc1}}]}}])-[d(x^{2_{[r_1]}}),x^{-\si_{r_1},1_{[\ol
r_{_{\ssc1}}]}}]
\vs{4pt}\\&
=2d(x^{2_{[r_1]}})-[d(x^{2_{[r_1]}}),x^{-\si_{r_1},1_{[\ol
r_{_{\ssc1}}]}}].
\end{array}
\eqno(3.20)$$%
Note from (3.18) that a nonzero term $x^{\g,\v k}$ appearing in
the left-hand side of (3.19) must satisfy $k_{\ol r_1}=\g_{r_1}$.
But such a term cannot appear in any image of $\ad
x^{-\si_{r_1},1_{[\ol r_{_{\ssc1}}]}}$, i.e., the coefficient of
such a term in the right-hand side of (3.19) is zero. Thus (3.19)
is in fact zero. Similarly, a nonzero term appearing in the
left-hand side of (3.20) must have the form
$x^{\a+2_{[r_1]}+\si_{r_1},\v i-1_{[\ol
r_{_{\ssc1}}]}}=x^{\a+1_{[r_1]},\v i-1_{[\ol r_{_{\ssc1}}]}}$ such
that $x^{\a,\v i}$ appears in $u$. But, if the term
$x^{\a+1_{[r_1]},\v i-1_{[\ol r_{_{\ssc1}}]}}$ appears in the
right-hand side of (3.20), then it must appear in
$d(x^{2_{[r_1]}})$, and by calculation, the coefficient of such a
term in the right-hand side of (3.20) is a multiple of%
$$2-((\a_{r_1}+1)-(i_{\ol r_1}-1))=i_{\ol r_1}-\a_{r_1}=0,%
$$%
i.e., (3.20) is also zero. Then (3.18)-(3.20) show that $u$
commutes with all elements in the set $\{x^{\b,\v
j},x^{2_{[r_1]}}\,|\,\b_{r_1}=j_{\ol r_1}\}$. Such an element $u$
has to be zero. This completes the proof of the claim.
\par%
{\bf Claim 4}.
Let $$A'_3=\{x^{2_{[p]}},t^{2_{[q]}}\,|\,p\in I_{1,5},q\in
I_6\},\eqno(3.21)$$%
and let $A_3=A'_3\cup\{x^{2_{[0]}}\}$ if $\G_0\ne\{0\}$ or
$A_3=A'_3\cup\{t_0\}$ if $\G_0=\{0\}$. We can suppose $d(A_3)=0$.
\par
Consider $d(x^{2_{[r]}})$ for $r\in I_1$. Note that $x^{2_{[r]}}$
is a common eigenvector for $\ad(A_1\cup A_2)$, and $d$ commutes
with $\ad(A_1\cup A_2)$. Thus $d(x^{2_{[r]}})$ is a common
eigenvector for $\ad(A_1\cup A_2)$. {}From this and (3.4)-(3.7), we
obtain that if $x^{\a,\v i}$ appears in $d(x^{2_{[r]}})$ with
coefficient $c_{\a,\v i}\ne0$, then
$$
\begin{array}{l}
\a_0=i_0=2+\a_{\ol r}-\a_r=0,\vs{4pt}\\
\a_{\ol p_1}-\a_{p_1}=i_{\ol p_2}=i_{p_3}=i_{\ol
p_4}-i_{p_4}=0,\\ \end{array}%
\eqno(3.22)$$
for $p_1\in I_{1,3}\bs\{r\},p_2\in I_{2,3},p_3\in I_3\cup
I_5,p_4\in I_6.$ If $\a_{\ol r}\ne-1$, we take $v=(\a_{\ol
r}+1)^{-1}c_{\a,\v i}x^{\a-2_{[r]}-\si_{[r]}}$ and replace $d$ by
$d+\ad v$. Since $v$ commutes with $A_1\cup A_2$ by (3.22), after
this replacement, we still have Claims 2 and 3, but then the term
$x^{\a,\v i}$ disappears in $d(x^{2_{[r]}})$. Thus we can suppose
$$
\a_{\ol r}=-1\mbox{ \ if \ }x^{\a,\v i}\mbox{ \ appears \,in \ }d(x^{2_{[r]}}).%
\eqno(3.23)$$%
Noting that $A_3$ is commutative, we have
$$
[y,d(x^{2_{[r]}})]=-[d(y),x^{2_{[r]}}]\;\;\;\;\for\;\;\;\;y\in A_3.%
\eqno(3.24)$$
Note from (3.21) and (3.23) that any nonzero term $x^{\b,\v j}$
appearing in the left-hand side of (3.24) must satisfies $\b_{\ol
r}=-1$. But such  a term does not appear in any image of $\ad
x^{2_{[r]}}$, thus (3.24) is in fact zero. Hence $d(x^{2_{[r]}})$ is
a common eigenvector for $\ad(A_1\cup A_2\cup A_3)$. {}From this and
formula (3.2), one can easily deduce that $d(x^{2_{[r]}})=0$.
Similarly, we can prove $d(y)=0$ for any $y\in A_3$. Thus the
claim follows.
\par
{\bf Claim 5}. Let
$$A'_4=\{x^{2_{[\ol p]}},t^{2_{[\ol q]}}\,|\,p\in I_{1,3},q\in
I_{4,6}\},
\eqno(3.25)$$
and let $A_4=A'_4\cup\{x^{-2_{[0]}}\}$ if $\G_0\ne\{0\}$ or
$A_4=A'_4$ if $\G_0=\{0\}$. We have $d(A_4)=0$.
\par
Set $B_3=A_1\cup A_2\cup A_3$. Consider $d(x^{2_{[\ol p]}}),\,p\in
I_{1,3}$. Note from formula (3.2) that $x^{2_{[\ol p]}}$ commutes
with all elements of $B_3$ except two elements, i.e.,
$$
[y,x^{2_{[\ol p]}}]=0,\;\;\; [x^{-\si_p},x^{2_{[\ol
p]}}]=2x^{2_{[\ol p]}},\;\;\; [x^{2_{[p]}},x^{2_{[\ol
p]}}]=4x^{-\si_p},
\eqno(3.26)$$
where $y\in
A'_3\bs\{x^{-\si_p},x^{2_{[p]}}\}$
(cf.~(2.7)). Applying $d$ to (3.26), we obtain that $d(x^{2_{[\ol
p]}})$, if not zero, is a common eigenvector for $\ad B_3$. But
from (3.2), we see that $\ad B_3$ does not have a common
eigenvector. Thus $d(x^{2_{[\ol p]}})=0$. Similarly we have
$d(y)=0$ for any $y\in A_4$. This proves the claim.
\par
{\bf Claim 6}. For any $p\in I_{1,3}$, let
$\G'_p=\{k_{[p]}+(2-k)_{[\ol p]}\,|\,k\in\Z\}\subset\G$. We can
suppose $d(x^\a)=0$ for $\a=k_{[p]}+(2-k)_{[\ol p]}\in\G'_p$
(recall notation (2.13)).
\par
Note that if $k=0,1,2$, then $\a=2_{[\ol p]},-\si_p,2_{[p]}$
respectively, and the result follows from Claims 2-5. Now the
general result can be proved as in the proof of Claim 5 by
induction on $|a|$.
\par
{\bf Claim 7}. We can suppose $d=0$.
\par
For any $\v i\in\JJ$, we define%
$$|\v i|=\sum_{p\in \wh J\bs J_6}i_p\mbox{ \ \ \ if \ \ \
}\G_0\ne\{0\},\mbox{ \ \ \ or \ \ \ }|\v i|=\sum_{p\in J\bs
J_6}i_p\mbox{
\ \ \ if \ \ \ }\G_0=\{0\}.%
$$%
For $n\in\N$, we denote
$$\KK_n={\rm span}\{x^{\a,\v i}\,\;|\;\,|\v i|\le n\}.%
\eqno(3.27)$$%
Inductively suppose we have proved $d(\KK_{n-1})=0$. Consider
$u=d(x^{\a,\v i})$ with $|\v i|=n$ and suppose a nonzero term
$x^{\b,\v j}$ appears in $u$.
\par
Assume that $p\in I_6$. Applying $d$ to
$$
[t^{1_{[p]}+1_{[p]}},x^{\a,\v i}]=(i_{\ol p}-i_p)x^{\a,\v
i},\;\;\; [t^{2_{[\ol p]}},[t^{2_{[p]}},x^{\a,\v
i}]]=-4(i_p+1)i_{\ol p}x^{\a,\v i}, \eqno(3.28)$$
we obtain
$$
[t^{1_{[p]}+1_{[p]}},u]=(i_{\ol p}-i_p)u,\;\;\; [t^{2_{[\ol
p]}},[t^{2_{[p]}},u]]=-4(i_p+1)i_{\ol p}u.%
\eqno(3.29)$$
By computing the coefficients of $x^{\b,\v j}$ in both sides of
each equation in (3.29), we obtain%
$$j_{\ol p}-j_p=i_{\ol p}-i_p,\;\;\;\;(i_p+1)i_{\ol
p}=(j_p+1)j_{\ol p},%
$$%
from this, since $i_i,i_{\ol p},j_p,j_{\ol p}$ are all nonnegative
integers, we obtain $$i_p=j_p,i_{\ol p}=j_{\ol p}.\eqno(3.30)$$
Similarly, we can prove%
$$
\b_0=\a_0,\;j_0=0 \mbox{ \ if \ }\G_0\ne\{0\},\mbox{ \ or \
}j_0=i_0\mbox{ \ if \ }\G_0=\{0\},\vs{-3pt}\eqno(3.31)$$
$$
\b_{\ol p}-\b_p=\a_{\ol p}-\a_p,\;(\b_p+1)\b_{\ol
p}=(\a_p+1)\a_{\ol p},\;j_p=j_{\ol p}=0\mbox{  for }p\in I_{1,3},\vs{-3pt}\eqno(3.32)$$
$$
\b_p=\a_p,\;\;\;j_p=j_{\ol p}=i_{\ol p}=0\mbox{ \ \ \ for \ }p\in
I_{4,5},\eqno(3.33)$$
and similarly, we have %
$$
\begin{array}{ll}
&(-(\b_{\ol p}-2)-3(\b_p+2))(\b_{\ol p}-1)\b_{\ol p}
\vs{4pt}\\
=
\!\!\!\!&
(-(\a_{\ol
p}-2)-3(\a_p+2))(\a_{\ol p}-1)\a_{\ol p}
\mbox{ \ \ for \ $p\in I_{1,3},$}
\end{array}
\eqno(3.34)
$$
which is the eigenvalue for $(\ad x^{-1_{[p]}+3_{[\ol p]}})(\ad
x^{2_{[p]}})^2$ corresponding to the eigenvector $u$, and
$$
\a_p=\a_{\ol p}=0\ \ \Rar\ \ \b_p=\b_{\ol
p}=0\;\;\;\;\for\;\;\;\;p\in I_{1,3},%
\eqno(3.35)$$%
because if $\a_p=\a_{\ol p}=0$ then $[x^{2_{[p]}},x^{\a,\v
i}],[x^{2_{[\ol p]}},x^{\a,\v i}]\in\KK_{n-1}$ and so $u$ commutes
with both $x^{2_{[p]}}$ and $x^{2_{[\ol p]}}$. {}From (3.32), we
solve
that%
$$(\b_p,\b_{\ol p})=(\a_p,\a_{\ol p})\mbox{ \  or \ }(-\a_{\ol
p}-1,-\a_p-1).%
\eqno(3.36)$$%
If $(\b_p,\b_{\ol p})\ne(\a_p,\a_{\ol p})$, by (3.34) and (3.36),
we obtain
$$(\a_p+\a_{\ol p})(\a_p+\a_{\ol
p}+1)(\a_p+\a_{\ol p}+2)=0\mbox{ \ if \ }(\b_p,\b_{\ol p})=(-\a_{\ol p}-1,-\a_p-1).%
\eqno(3.37)$$
%
%
%
Thus if we suppose%
$$
\a_q+\a_{\ol q}\ne0,-1,-2\mbox{ \ \ for \,all \ \ }q\in \{p\in I_{1,3}\,|\,(\a_p,\a_{\ol p})\ne0\},%
\eqno(3.38)$$%
then we obtain from (3.35)-(3.37) that $\b_p=\a_p$ for $p\in
J_{1,3}$. This together with (3.30)-(3.33) gives that $(\b,\v
j)=(\a,\v i_{J_6})$ if $\G_0\ne\{0\}$ or $(\b,\v j)=(\a,\v i_{\wh
J_6})$ if $\G_0=\{0\}$. i.e., $u$ has at most one term, so we can
suppose%
$$
d(x^{\a,\v i})=c_{\a,\v i\,}x^{\a,\v i_{J_6}}\mbox{ if
}\G_0\ne\{0\},\mbox{ \ or \ }d(x^{\a,\v i})=c_{\a,\v i\,}x^{\a,\v
i_{\wh J_6}}\mbox{ if
 }\G_0=\{0\},%
\eqno(3.39)$$%
for some $c_{\a,\v i}\in\F$. Note that we deduce (3.39) under the
condition (3.38).
For $\a\in\G$, let%
$$S_\a=\{p\in I_{1,3}\,|\,(\a_p,\a_{\ol p})\ne0,\;
\a_p+\a_{\ol p}=0,-1,-2\}\mbox{ \ and \ }
m_\a=|S_\a|.%
\eqno(3.40)$$%
We want to prove by induction on $m_\a$ the following
\par
{\bf Statement 1}. Condition (3.38) can be removed, i.e., (3.39)
holds for all $\a\in\G$.
\par%
Note that if we can find $(\mu,\v j),(\nu,\v k)\in\G\times\JJ,\,|\v j|\le n,|\v k|\le n$ such that%
$$[x^{\mu,\v j},x^{\nu,\v
k}]=c_1x^{\a,\v i}+c_2x^{\a-\si_p,\v i}+y,%
\eqno(3.41)$$ %
for some $c_1,c_2\in\F,\,p\in S_\a,y\in\KK_{n-1}$ with
$c_1\ne0,m_\mu=0,m_\nu<m_\a$, then we are done because either
$p\notin S_{\a-\si_p}$ (then $m_{\a-\si_p}=m_\a-1$), or $p\in
S_{\a-\si_p}$ (then $\a_p+\a_{\ol p}=-2$ and
$(\a-\si_p)_p+(\a-\si_p)_{\ol p}=0$, and we can use induction on
$\a_p+\a_{\ol p}$), and we can use (3.1), the induction on $m_\a$
and the assumption that $d(\KK_{n-1})=0$ to obtain (3.39).
\par
Assume that $p\in S_\a$. In all cases below, we can choose $\v
j=0,\v k=\v i$ and choose $\mu,\nu$ as follows so that (3.41) can
be satisfied (where $j,k,l\in\Z$ are suitable integers to
guarantee $c_1\ne0,m_\mu=0,m_\nu<m_\a$),
$$
\mu=j_{[p]}+k_{[\ol p]},\;\;\nu=\a+(1-j)_{[p]}+(1-k)_{[\ol
p]}\mbox{ \ \ if
\ \ }(\a_p,\a_{\ol p})\ne(-1,-1),%
\eqno(3.42)$$
$$
\mu=k\si_p,\;\nu=\a-k\si_p\mbox{ \ \ \ if \ \ \
}\a_0\ne0,\;(\a_p,\a_{\ol
p})=(-1,-1),%
\eqno(3.43)$$
$$
\left\{\begin{array}{ll} \mu=\!\!\!\!&j_{[q]}+k_{[\ol
q]}+l\si_p,\;\nu=\a+(1-j)_{[q]}+(1-k)_{\ol [q]}-l\si_p
\vs{4pt}\\&
\mbox{if }\a_0=0,\;(\a_p,\a_{\ol p})=(-1,-1)
\vs{4pt}\\&
\mbox{and
}\exists\,q\in I_{1,3}\bs\{p\}\mbox{
with }(\a_q,\a_{\ol q})\ne(-1,-1),\end{array}%
\right.\eqno(3.44)$$
$$
\left\{\begin{array}{ll}
\mu=\!\!\!\!&2_{[0]},\;\nu=\a-2_{[0]}
\vs{4pt}\\&
\mbox{if
}\G_0\ne\{0\},\a_0=0,\vt(\a,\v
i)\ne4,\;\a_q=-1\mbox{ for all }q\in J_{1,3}.%
\end{array}\right.
\eqno(3.45)$$%
Now we are left the following two cases:
$$
\G_0\ne\{0\},\a_0=0,\vt(\a,\v i)=4,\;\a_q=-1\mbox{ for all }q\in J_{1,3}.%
\eqno(3.46)$$
$$
\G_0=\{0\},\;\a_q=-1\mbox{ for all }q\in J_{1,3}.%
\eqno(3.47)$$
Suppose we are in case (3.46). If $x^\b$ appears in $u$ with
nonzero coefficient, and $(\b_p,\b_{\ol p})\ne(\a_p,\a_{\ol p})$
for some $p\in I_{1,3}$, then (3.36) and (3.46) show that
$(\b_p,\b_{\ol p})=(0,0)$, while $(\a_p,\a_{\ol p})=(-1,-1)$. Let
$m$ be the number of such $p$'s, then (3.30) and (3.33) show that
$$\vt(\b,0)=\vt(\a,0)+2m=\vt(\a,\v i)+2m.
\eqno(3.48)$$ %
(Recall definition of $\vt(\a,\v i)$ in (3.3).) But $u$ is also an
eigenvector for \linebreak $(\ad x^{2_{[0]}})(\ad x^{-2_{[0]}})$ (cf.~Claims
4 and 5) with eigenvalue equal to
$$
-4(2-\vt(\a,\v i))(4-\vt(\a,\v i))=-4(2-\vt(\b,0))(4-\vt(\b,0)),%
\eqno(3.49)
$$%
(cf.~(3.2)), from this we obtain $\vt(\a,\v i)=\vt(\b,0)$ or
$\vt(\a,\v i)+\vt(\b,0)=6$, and we get a contradiction from this
and (3.46), (3.48). Thus we must have $\b=\a$, i.e., (3.39)
holds in case (3.46). Now suppose we are in case (3.47). Then
we have $x^{\a,\v i}=c_1[x^{k\si_p,1_{[0]}},x^{\a-k\si_p,\v i}]$
for some $c_1\ne0,k\in\Z$ so that (3.41) is satisfied. This
completes the proof of Statement 1.
\par%
Now we prove $d(x^{\a,\v i})=0$ by considering the following cases.%
\par%
{\it Case 1}: $\G_0\ne\{0\}$ and $\v i=\v i_{J_6}$ (i.e., $n=0$).
\par%
We want to prove
$$
c_{\a+\b,\v i+\v j}=c_{\a,\v i}+c_{\b,\v
j}\;\;\;\;\for\;\;\;\;\a,\b\in\G\mbox{ \ and \ }\v i=\v i_{J_6},\v
j=\v j_{J_6}\in J,
\eqno(3.50)$$
%
We define
$$\phi_0(\a,\v i,\b,\v j)=(2-\vt(\a,\v i))\b_0-\a_0(2-\vt(\b,\v j)).%
\eqno(3.51)$$%
Applying $d$ to (3.2), and comparing the coefficients of
$x^{\a+\b,\v i+\v j}$ in both sides, we obtain that (3.50) holds
if $\phi_0(\a,\v i,\b,\v j)\ne0$. Assume that $\phi_0(\a,\v
i,\b,\v j)=0$. If $\a_0,\b_0\ne0$, then we can always choose
$\g\in\G$, such that
$$
\phi_0(\a+\g,\v i,\b-\g,\v j),\ \phi_0(\a,\v i,\g,0),\
\phi_0(\b,\v j,-\g,0),\
\phi_0(\g,0,-\g,0)\ne0,%
\eqno(3.52)$$%
so we have (note that $c_{0,0}=0$ since $d(1)=0$)
$$\begin{array}{ll}
c_{\a+\b,\v i+\v j}\!\!\!\!&=c_{(\a+\g)+(\b-\g),\v i+\v
j}=c_{\a+\g,\v i}+c_{\b-\g,\v j}\vs{4pt}\\
&=c_{\a,\v i}+c_{\b,\v j}+c_{\g,0}+c_{-\g,0}
\vs{4pt}\\&
=c_{\a,\v i}+c_{\b,\v
j}+c_{\g+(-\g),0}=c_{\a,\v i}+c_{\b,\v j},\end{array}
\eqno(3.53)$$
i.e., (3.50) holds if $\a_0,\b_0\ne0$. If $\a_0\ne0=\b_0$, we can
choose $\g\in\G$ such that $\a_0\ne\g_0\ne0$, so
$$\begin{array}{ll}
c_{\a,\v i}+c_{\b,\v j}\!\!\!\!&=c_{\a,\v i}+c_{(\b-\g)+\g,\v
j}=c_{\a,\v i}+c_{\b-\g,\v j}+c_{\g,0}\vs{4pt}\\
&=c_{\a+\b-\g,\v i+\v j}+c_{\g,0}=c_{\a+\b,\v i+\v
j}+c_{-\g,0}+c_{\g,0}=c_{\a+\b,\v i+\v j},\end{array}
$$
i.e., (3.50) holds if $\a_0\ne0=\b_0$. If $\a_0=\b_0=0$, choose
$\g\in\G$ such that $\g_0\ne0$, then we have all equalities of
(3.53). Hence (3.50) holds in this case. Now (3.50) gives %
$$c_{\a,\v i}=c_{\a,0}+\sum_{p\in J_6}i_pc_{0,1_{[p]}}
\;\;\;\;\for\;\;\;\;\a\in\G\mbox{ \ and \ }\v i=\v i_{J_6}\in J,\eqno(3.54)$$%
But (3.53), (3.21) and (3.25) show that
$2c_{0,1_{[p]}}=c_{0,2_{[p]}}=0$ for $p\in J_6$. Thus $c_{\a,\v
i}=c_{\a,0}$ and (3.50) shows that $\mu:\a\mapsto c_{\a,0}$
defines an element $\mu\in\HOM$, and by substituting $d$ by
$d-d_\mu$, we can suppose $c_{\a,0}=0$ for all $\a\in\G$, i.e., we
proved that $d(x^{\a,\v i})=0$ for all $\a\in\G,\,\v i=\v
i_{J_6}\in J$.
\par
{\it Case 2}: $\G_0=\{0\},\v i=\v i_{\wh J_6}$ (i.e., $n=0$).%
\par
Applying $d$ to $[1,x^{\a,\v i}]=2i_0x^{\a,\v i-1_{[0]}}$ we
obtain by induction on $i_0$ that $c_{\a,\v i}=c_{\a,\v i_{J_6}}$.
Applying $d$ to (3.2) and comparing coefficient of $x^{\a+\b,\v
i+\v j-1_{[0]}}$, we obtain $c_{\a,\v i}+c_{\b,\v j}=c_{\a+\b,\v
i+\v j}$ if $(2-\vt(\a,\v i))j_0-i_0(2-\vt(\b,\v j))\ne0$. But
$i_0,j_0\in\N$ are arbitrary, we have $c_{\a,\v i}+c_{\b,\v
j}=c_{\a+\b,\v i+\v j}$ if $(\vt(\a,\v i),\vt(\b,\v j))\ne(2,2)$.
The rest of the proof is the same
as in Case 1.%
\par%
{\it Case 3}: $\G_0\ne\{0\}$ and $n>0$.
\par
%
Note that a general element $x^{a\,\v i}$ can be generated by the set%
$$
S=\{x^{\b,\v j},t_p\,|\,\b\in\G,\,\v j=\v j_{J_6},\;\;p\in\ol I_2\cup J_3\cup I_5\}.%
\eqno(3.55)$$%
Recall (3.39). By replacing $d$ by $d-\sum_{p\in\ol I_2\cup
J_3\cup I_5}c_{{\ssc\,}0,1_{[p]}}\ptl_{t_p}$, we can suppose
$d(t_p)=0$ for $p\in\ol I_2\cup J_3\cup I_5$ (this replacement
does not affect the fact that $d(x^{\a,\v i})=0$ for $a\in\G,\,\v
i=\v i_{J_6}$). Thus $d(S)=0$ and so $d(x^{\a,\v
i})=0$.\par%
{\it Case 4}: $\G_0=\{0\}$ and $n>0$.%
\par
Similar to Case 3, 
a general element $x^{\a,\v i}$ 
can be generated by $S'=\{x^{\b,\v j},t_p\,|\,\b\in\G,\,\v j=\v
j_{\wh J_6},\;\;p\in\ol I_2\cup J_3\cup I_5\}$. The rest of the
proof of this case is exactly the
same as that of Case 3.%
\par%
This proves
Claim 7 and Theorem 3.1.%
\qed\par\
\vs{-2pt}\par%
\cl{\bf4. \ Second cohomology groups}%
\par
In this section, we shall determine the second cohomology
groups of the contact Lie algebra $\KK=\KK(\v\ell,\G)$. Recall
that a {\it 2-cocycle} on $\KK$ is an $\mbb{F}$-bilinear function
$\psi:\KK\times \KK\rar \mbb{F}$ satisfying the following
conditions:
$$
\psi(v_1,v_2)=-\psi(v_2,v_1)\mbox{\ \ (skew-symmetry)},
\eqno(4.1)$$
$$
\psi([v_1,v_2],v_3)+\psi([v_2,v_3],v_1)+\psi([v_3,v_1],v_2)
=0
\mbox{\ \ (Jacobian identity)},
\eqno(4.2)$$
for $v_1,v_2,v_3\in \KK$.
Denote by $C^2(\KK,\mbb{F})$ the vector space of 2-cocycles on $\KK$.
For any $\mbb{F}$-linear function $f:\KK\rar\mbb{F}$, one can define a 2-cocycle
$\psi_f$ as follows
$$
\psi_f(v_1,v_2)=f([v_1,v_2])\;\;\;\for\;\;\;v_1,v_2\in \KK.
\eqno(4.3)$$
Such a 2-cocycle is called a {\it 2-coboundary} or a {\it trivial
2-cocycle} on $\KK$.
Denote by $B^2(\KK,\mbb{F})$ the
vector space of 2-coboundaries on $\KK$.
A 2-cocycle $\phi$ is said to be
{\it equivalent to} a 2-cocycle $\psi$ if $\phi-\psi$ is trivial. For
a 2-cocycle $\psi$, we denote by $[\psi]$ the equivalent class of $\psi$.
The quotient space
\\[4pt]\hs{1ex}$
H^2(\KK,\mbb{F})\!=\!C^2(\KK,\mbb{F})/B^2(\KK,\mbb{F})\!
=\{\mbox{the equivalent classes of 2-cocycles}\},
$\hfill(4.4)\\[6pt]
is called the {\it second cohomology group} of $\KK$.
\par\ni
{\bf Lemma 4.1}. {\it If $\JJ_0=\N$ or $\iota_6\ne\ell_1$, then
$H^2(\KK,\mbb{F})=0$.}
\par\ni
{\it Proof.} Let $\psi$ be a 2-cocycle. Say $\ell_2\ne0$ (from the
proof below, one sees that the proof is exactly similar if
$\JJ_0=\N$ or $\ell_i\ne0$ for some $i\ne1,2$). Fix $p\in I_2$ and
define a linear function $f$ by induction on $i_{\ol p}$ as
follows:
$$
f(x^{\a,\v i})=\left\{
\begin{array}{lll}
(\a_{\ol
p}-\a)^{-1}(\psi(x^{-\si_p},x^{\a,\v i}) -i_{\ol p}f(x^{\a,\v
i-1_{[\ol p]}}))\hfill&\mbox{if \ }\a_{\ol p}\ne\a_p,
\vs{4pt}\hfill\cr (i_{\ol p}+1)^{-1}\psi(x^{-\si_p},x^{\a,\v
i+1_{[\ol p]}})\hfill& \mbox{if \ }\a_{\ol p}=\a_p,
\end{array}\right.
\eqno(4.5)$$ for $(\a,\v i)\in\G\times\JJ$. Set
$\phi=\psi-\psi_f$. Then (3.5), (4.3) and (4.5) show that
$$
\phi(x^{-\si_p},x^{\a,\v i})=0\;\;\;\for\;\;\;(\a,\v
i)\in\G\times\JJ. \eqno(4.6)$$ Using Jacobian identity (4.2), we
obtain
$$
\begin{array}{lll}
0\!\!\!\!&=\phi(x^{-\si_p},[x^{\a,\v i},x^{\b,\v
j}])\hfill\vs{4pt}\cr&
= (\a_{\ol p}+\b_{\ol p}
-\a_p-\b_p)
\phi(x^{\a,\v i},x^{\b,\v j})
\vs{4pt}\\&\ \ \
+i_{\ol p}\phi(x^{\a,\v
i-1_{[\ol p]}},x^{\b,\v j})
 +
j_{\ol p}\phi(x^{\a,\v i},x^{\b,\v
j-1_{[\ol p]}}),
\end{array}
\eqno(4.7)$$
for $(\a,\v i),\,(\b,\v
j)\in\G\times\JJ$. If $\a_{\ol p}+\b_{\ol p}\ne\a_p+\b_p$, by
induction on $i_{\ol p}+j_{\ol p}$, we obtain $\phi(x^{\a,\v
i},x^{\b,\v j})=0$. Otherwise (4.7) gives
$$
\phi(x^{\a,\v i},x^{\b,\v j})
=-j_{\ol p}(i_{\ol p}+1)^{-1}\phi(x^{\a,\v i+1_{[\ol p]}},
x^{\b,\v j-1_{[\ol p]}}),
\eqno(4.8)$$
and by induction on $j_{\ol p}$, we again have
$\phi(x^{\a,\v i},x^{\b,\v j})=0$. Thus $\phi=0$.
\qed\par
{}From now on, we assume that $\JJ_0=\{0\}$ and $\iota_6=\ell_1$.
Denote
$$\si=\sum_{p\in
I_1}\si_p=(0,-1,-1,...,-1)\in\F^{1+2\ell_1},%
\eqno(4.9)$$
as in (2.19). Suppose $\psi$ is a 2-cocycle. We define a linear
function $f$ as follows: First for $\a\in\G\bs\{\si\}$ with
$\a_0=0$, we define
$$
p_\a={\rm min}\{p\in I_1\,|\,(\a_p,\a_{\ol p})\ne(-1,-1)\}.%
\eqno(4.10)
$$%
Then we set
$$f(x^\a)=\left\{\begin{array}{ll}
(4(2+\ell_1))^{-1}\psi(x^{-2_{[0]}},x^{2_{[0]}+\si})\!\!\!\!&\mbox{if \
}\a=\si,\vs{4pt}\\ %
(2\a_0)^{-1}\psi(1,x^\a)&\mbox{if \
}\a_0\ne0,\vs{4pt}\\ %
(\a_{\ol p}-\a_p)^{-1}\psi(x^{-\si_p},x^\a)&\mbox{if \
}\a\ne\si,\,\a_0=0,\a_p\ne\a_{\ol p},\vs{4pt}\\ %
-(2(\a_{\ol p}+1))^{-1}\psi(x^{-2_{[p]}},x^{\a-1_{[p]}+1_{[\ol
p]}})\!\!\!\!\!\!\!\!\!\!\!\!\!\!\!\!\!\!\!\!\!\!\!\!\!\!\!\!\!\!\!\!
\vs{2pt}\\&\hs{-10ex}\hfill\mbox{if \
}\a\ne\si,\,\a_0=0,\a_p=\a_{\ol p}\ne-1,%
\\ \end{array}\right.
\eqno(4.11)$$ where $p=p_\a$. Set $\phi=\psi-\psi_f$.%
\par\ni%
{\bf Lemma 4.2}. {\it For $a\in\G$, we have%
$$%
\begin{array}{llll}
\phi(1,x^\a)&=&0,\vs{4pt}\\ %
\phi(x^{-\si_p},x^\a)&=&0&\mbox{ \ for \ \ \ }p\in I_1,\vs{4pt}\\
\phi(x^{2_{[q]}},x^\a)&=&0&\mbox{ \ for \ \ \ }q\in J_1.%
\end{array}%
\eqno\begin{array}{r}\vs{4pt}(4.12)\!\\ \vs{4pt}(4.13)\!\\
(4.14)\!\end{array}
$$}%
\par\ni{\it Proof.} If $\a_0\ne0$, by definition, we have
$\phi(1,x^\a)=\psi(1,x^\a)-f([1,x^\a])=0$. If $\a_0=0$, since
$\KK=[\KK,\KK]$, we can always write $x^\a$ as a linear
combination of $[x^\b,x^\g]$ such that $\b_0+\g_0=0$. Using
Jacobian identity (4.2), we have
$\phi(1,[x^\b,x^\g])=(\b_0+\g_0)\phi(x^\b,x^\g)=0$. This proves
(4.12). For any $\a,\b\in\G$, from (4.12), we obtain
$0=\phi(1,[x^\a,x^\b])=(\a_0+\b_0)\phi(x^\a,x^\b)$, thus
$$
\phi(x^\a,x^\b)=0\;\;\;\;\for\;\;\;\;\a,\b\in\G\;\;\;\;\mbox{with}\;\;\;\a_0+\b_0\ne0.%
\eqno(4.15)%
$$%
Consider (4.13). We can suppose $\a_0=0$ by (4.15). If
$\a_p=\a_{\ol p}$, we have (4.13) by writing $x^\a$ as a linear
combination of $[x^\b,x^\g]$ with $\b_p+\g_p=\b_{\ol p}+\g_{\ol
p}$ and using (4.2). Assume $\a_p\ne\a_{\ol p}$. Then $p_a\le p$.
If $p=p_\a$, we have (4.13) by definition (4.11). If $r:=p_a<p$,
then either $x^\a=(\a_{\ol r}-\a_r)^{-1}[x^{-\si_r},x^\a]$ (if
$\a_{\ol r}\ne \a_r$) or $x^\a=(2(\a_{\ol
r}+1))^{-1}[x^{2_{[r]}},x^{\a-1_{[r]}+1_{[\ol r]}}]$ (if
$\a_r=\a_{\ol r}\ne-1$), and we have (4.13) by Jacobian identity
(4.2) and definition of $f(x^\a)$ in (4.11). This proves (4.13).
As in (4.15), we have
$$
\phi(x^\a,x^\b)=0\;\;\;\;\mbox{if}\;\;\;a_p+\b_p\ne\a_{\ol
p}+\b_{\ol p}\;\;\;\mbox{for \ some}\;\;\;p\in I_1.%
\eqno(4.16)$$%
Consider (4.14). First suppose $q\in I_1$. If $\a_{\ol q}=0$, then
as above by writing $x^\a$ as a sum of the form $[x^\b,x^\g]$ with
$\b_{\ol q}=\g_{\ol q}=0$, we have (4.14). Assume $\a_{\ol
q}\ne0$. Then $p_{\a+1_{[q]}-1_{[\ol q]}}\le q$. By (4.15)-(4.16)
and the definition of $f(x^{\a+1_{[q]}-1_{[\ol q]}})$, using the
same arguments as above, we have (4.14). Finally suppose $q=\ol
p\in\ol I_1$. We can assume that $\a_0=0,\a_r=\a_{\ol r}=-1$ for
$r\in I_1\bs\{p\}$ and $\a_p=\a_{\ol p}+2$ and $\a_p=-1$
(otherwise either we can use (4.15)-(4.16) to obtain (4.14) or we
can write $x^\a$ as $[x^{2_{[r]}},y]$ for some $r\in I_1,y\in\KK$
and use (4.2), (4.13) to obtain (4.14)). But then $x^\a$ is a
multiple of $[x^{-2_{[0]}},x^{2_{[0]}+\a}]$ and $[x^{2_{[\ol
p]}},x^{-2_{[0]}+\a}]$ is a multiple of $x^{-2_{[0]}+\si}$, and we
can use (4.2) and the first case of (4.11) to obtain (4.14). This
proves the lemma.
\par\ni%
{\bf Lemma 4.3}. {\it We have $\phi(x^{i_{[p]}+(2-i)_{[\ol
p]}},x^\a)=0$ for all $\a\in\G,i\in\Z,p\in I_1$.}
\par\ni%
{\it Proof.} We have the result for $i=0,1,2$ by (4.13)-(4.14).
The general result can be obtained by writing
$x^{i_{[p]}+(2-i)_{[\ol p]}}$ as a multiple of the form
$[x^{2_{[r]}},x^{(k_{[p]}+(2-k)_{[\ol p]}}]$ for some $r=p,\ol p$
and some $k\in\Z$ such that $|k|<|i|$, and using (4.2) and by
induction on
$|i|$.%
\qed\par\ni%
{\bf Lemma 4.4}. {\it We have
$\phi(x^\a,x^\b)=\d_{\a+\b,\si\,}c_\a$
for all $\a,\b\in\G$ and some $c_\a\in\F$.}%
\par\ni%
{\it Proof.} Assume that $\phi(x^\a,x^\b)\ne0$. First we prove%
$$
\a_p=\a_{\ol p}=0 \ \ \Rar \ \ \b_p=\b_{\ol
p}=-1\;\;\;\;\for\;\;\;\;p\in I_1.%
\eqno(4.17)
$$%
Suppose $\a_{\ol p}=0,\b_{\ol p}\ne-1$, then
$[x^{2_{[p]}},x^\a]=0$, we have
$$
\begin{array}{ll}
0
\!\!\!\!&
=
\phi([x^{2_{[p]}},x^\a],x^{\b-1_{[p]}+1_{[\ol
p]}})
\vs{4pt}\\&
=-\phi(x^\a,[x^{2_{[p]}},x^{\b-1_{[p]}+1_{[\ol p]}}])
=-(\b_{\ol p}+1)\phi(x^\a,x^\b),
\end{array}
\eqno(4.18)
$$
which is a contradiction with the assumption that
$\phi(x^\a,x^\b)\ne0$. Thus $\a_{\ol p}=0$ implies $\b_{\ol
p}=-1$. Similarly $\a_p=0$ implies $\b_p=-1$, i.e., we have
(4.17). Now set
$$S=\{x^\a\,|\,\a\in\G,\,(\a_p,\a_{\ol p})\ne0,\;\a_p+\a_{\ol p}\ne0,-1,-2\mbox{ \ for \,all \ }p\in I_1\},%
\eqno(4.19)$$%
 First assume that $x^a\in S$. By (4.15)-(4.16), we can suppose
$$
\a_0+\b_0=0,\;\;\;\a_p+\b_p=\a_{\ol p}+\b_{\ol
p}\;\;\;\;\for\;\;\;\;p\in I_1.%
\eqno(4.20)$$
By Lemma 4.3 and Jacobian identity (4.2), we have%
$$
\begin{array}{ll}
-4(\a_p+1)\a_{\ol p}\phi(x^\a,x^\b)\!\!\!\!&=\phi([x^{2_{[\ol
p]}},[x^{2_{[p]}},x^\a]],x^\b)\vs{4pt}\\
&=\phi(x^\a,[x^{2_{[p]}},[x^{2_{[\ol p]}},x^\b]])
\vs{4pt}\\&
=-4(\b_{\ol
p}+1)\b_p\phi(x^\a,x^\b),
\end{array}%
\eqno(4.21)$$
and%
$$
\begin{array}{ll}
&4(-(\a_p-2)-3(\a_p+2))(\a_{\ol p}-1)\a_{\ol
p}\phi(x^\a,x^\b)\vs{4pt}\\ =&\phi([x^{-1_{[
p]}+3_{[\ol p]}},[x^{2_{[p]}},[x^{2_{[p]}},x^\a]]],x^\b)\vs{4pt}\\
=&-\phi(x^\a,[x^{2_{[p]}},[x^{2_{[p]}},[x^{-1_{[p]}+3_{[\ol p]}},x^\b]]])\vs{4pt}\\
=&-4(\b_{\ol p}+1)(\b_{\ol p}+2)(-\b_{\ol
p}-3\b_p)\phi(x^\a,x^\b),
\end{array}%
\eqno(4.22)$$
for $p\in I_1$. If $\phi(x^\a,x^\b)\ne0$, from (4.17),
(4.20)-(4.22) and the fact that $\a\in S$, we obtain
$\a_p+\b_p=-1$ for all $p\in J_1$, i.e., $\a+\b=\si$. Thus we have
the result of Lemma 3.4 in the case $\a\in S$. The general result
of Lemma 4.4 follows from Jacobian identity (4.2) and the fact
that $S$ is a generating set of $\KK$ (the proof of Statement 1 in
Section 3 shows
 that $S$ is a generating set of $\KK$).%
\qed\par\ni%
{\bf Lemma 4.5}. {\it We have $c_\a=0$ for
$\a\in\G$.}%
\par\ni%
{\it Proof.} Fix $p\in I_1$ (note that we require at the beginning
of Section 2 that $\iota_6>0$, i.e., $\ell_1>0$ by the assumption
$\iota_6=\ell_1$) and take $\b=1_{[0]}+1_{[\ol p]}$. For any
$\a\in\G$ with $\a_p\ne0$, considering
$\phi([x^\a,x^\b],x^{-\a-\b+\si-\si_p})$ and using (4.2), we have
$$
(\a_p\b_{\ol p}-\a_{\ol p}\b_p)c_{\a+\b+\si_p}=(\a_p\b_{\ol
p}-\a_{\ol p}\b_p)c_\a+(\a_p\b_{\ol p}-\a_{\ol p}\b_p)c_\b,
\eqno(4.23)
$$
that is
$$
c_{\a+\b+\si_p}=c_\a+c_\b.%
\eqno(4.24)$$
Now considering $\phi([x^\a,x^{\b+\si_p}],x^{-\a-\b+\si})$, we
have
$$\begin{array}{ll}
&
(2-\vt(\a,0)-3\a_0)c_{\a+\b+\si_p}
\vs{4pt}\\
&=
(-3(\a_0+1)-(1+\vt(\a,0)+2\ell_1))c_\a\vs{4pt}\\
&%
-((2-\vt(\a,0))(-\a_0-1)-\a_0(1+\vt(\a,0)+2\ell_1))c_{\b+\si_p}.
\end{array}%
\eqno(4.25)$$ Using (4.24) in (4.25), we obtain that%
$$%
2(3+\ell_1)c_\a=(2-\vt(\a,0)-3\a_0)c_\b-(\a_0(3+2\ell_1)+2-\vt(\a,0))c_{\b+\si_p}.\eqno(4.26)
$$%
Setting $\a=\si$ (then $\a_p\ne0$), by skew-symmetry (4.1) and the
second case of (4.11), we have $c_\si=-c_0=0$, we obtain
$c_\b=c_{\b+\si_p}$. Thus (4.26) becomes
$$c_\a=\a_0c_\b,%
\eqno(4.27)$$
holds for all $\a\in\G$ with $\a_p\ne0.$ If $\a_p=0$, then
$(-\a+\si)_p\ne0$ and by skew-symmetry (4.1), we have
$c_\a=-c_{-\a+\si_p}=-(-\a_0c_\b)=\a_0c_\b$, i.e., (4.27) holds
for all $\a\in\G$. But the first case of (4.11) and (4.27) show
that
$0=c_{-2_{[0]}}=-2c_\b$, and so (4.27) gives that $c_\a=0$ for all $\a\in\G$.%
\qed\par%
The above lemmas have proved the following theorem.
\par\ni%
{\bf Theorem 4.6}. {\it $H^2(\KK,\F)=0$.}%
\qed\par\ \vs{-10pt}\par%
\cl{\bf References} \small\par\ni%
\hi3.5ex\ha1 [1] R.~Farnsteiner, ``Derivations and central
extensions of finitely
  generated graded Lie algebras,'' {\it J.~Algebra} {\bf118} (1988), 33-45.
\par\ni\hi3.5ex\ha1
[2] V.~G.~Kac, ``A description of filtered Lie algebras whose associated
 graded Lie algebras are of Cartan types,'' {\it Math.~of USSR-Izvestijia}
 {\bf 8} (1974), 801-835.
\par\ni\hi3.5ex\ha1
[3] V.~G.~Kac, ``Lie superalgebras,'' {\it Adv.~Math.} {\bf 26} (1977), 8-96.
  \par\ni\hi3.5ex\ha1
[4] V.~G.~Kac, {\it Infinite Dimensional Lie Algebras}, 3rd ed.;
  Combridge Univ. Press, 1990.
 \par\ni\hi3.5ex\ha1
[5] J.~M.~Osborn, ``New simple infinite-dimensional Lie algebras of
 characteristic 0,'' {\it J.~Algebra} {\bf 185} (1996), 820-835.
 \par\ni\hi3.5ex\ha1
[6] J.~M.~Osborn, K.~Zhao, ``Generalized Cartan type $K$ Lie
algebras in
 characteristic 0,'' {\it Commun.~Algebra} {\bf 25} (1997), 3325-3360.
\par\ni\hi3.5ex\ha1
[7] Y.~Su, X.~Xu, ``Structure of contact Lie algebras related to
 locally-finite derivations,'' {\it Manuscripta Math.} {\bf112} (2003),
 231-257.
 \par\ni\hi3.5ex\ha1
[8] Y.~Su, X.~Xu and H.~Zhang, ``Derivation-simple algebras and the
 structures of Lie algebras of Witt type,'' {\it J.~Algebra} {\bf 233}
 (2000), 642-662.
 \par\ni\hi3.5ex\ha1
[9] X.~Xu, ``New generalized simple Lie algebras of Cartan type over a field
 with characteristic 0,''  {\it J.~Algebra} {\bf 224} (2000), 23-58.
 \par\ni\hi4ex\ha1
[10] X.~Xu, ``Quadratic conformal superalgebras,'' {\it J.~Algebra}
 {\bf 231} (2000), 1-38.
\end{document}